\documentclass[11pt]{article}
\usepackage{amsfonts}
\usepackage[margin=2cm]{geometry}
\usepackage{graphicx,color}
\usepackage[hang,small,bf]{caption}
\usepackage{xspace}
\usepackage{enumerate}
\usepackage{amssymb,amsmath}
\usepackage{amsthm}

\usepackage[T1]{fontenc}
\usepackage[utf8]{inputenc}
\usepackage{authblk}

\usepackage{natbib}
\bibpunct[, ]{(}{)}{,}{a}{}{,}%

\newtheorem{thm}{Theorem}

\newtheorem{lemma}[thm]{Lemma}

\newtheorem{definition}[thm]{Definition}

\begin{document}

	\title{Constant Factor Approximate Solutions  \\ for Expanding Search on General Networks}
	
	\author[]{Steve Alpern}
	\affil[]{Operations Research and Management Science Group, Warwick Business School, \\
			University of Warwick, Coventry CV4 7AL, United Kingdom, steve.alpern@wbs.ac.uk}
	\author[]{Thomas Lidbetter}
	\affil[]{Department of Management Science and Information Systems, Rutgers Business School, \\
			1 Washington Park, Newark, NJ 07102, USA, t.r.lidbetter@lse.ac.uk}

	\date{}
	\maketitle
	\thispagestyle{empty}
	\maketitle
	
	\begin{abstract}
	We study the classical problem introduced by R. Isaacs and S. Gal of
	minimizing the time to find a hidden point $H$ on a network $Q$ moving
	from a known starting point. Rather than adopting the traditional continuous
	unit speed path paradigm, we use the ``expanding search'' paradigm recently
	introduced by the authors. Here the regions $S\left( t\right) $ that have
	been searched by time $t$ are increasing from the starting point and have
	total length $t$. Roughly speaking the search follows a sequence of arcs
	$a_{i}$ such that each one starts at some point of an earlier one. This
	type of search is often carried out by real life search teams in the hunt for missing persons, escaped convicts, terrorists or lost airplanes. The paper
	which introduced this type of search solved the adversarial problem (where
	$H$ is hidden to take a long time to find) for the cases where $Q$ is a tree
	or is 2-arc-connected. This paper solves the game on some additional families of networks.
	However the main contribution is to give strategy classes which can be used
	on any network and have expected search times which are within a factor
	close to 1 of the value of the game (minimax search time). We identify cases
	where our strategies are in fact optimal. 
		
	\end{abstract}{\bf Keywords:} teams: games/group decisions; search/surveillance; tree algorithms; networks/graphs

\section{Introduction}

The problem of optimal search for a stationary hidden object $H$ (a 
\textit{Hider}) on a network with given arc lengths goes back to the early
work of \cite{Isaacs} and \cite{Gal79}. In the traditional version, the
Searcher proceeds continuously and at unit speed from a known starting point 
$O$ on the network $Q$ until reaching $H$ at some time $T$. The aim is find a randomized search that minimizes this time $T$,
the {\em search time}, in the worst case. An equivalent approach, which we prefer to adopt, is to study a zero-sum game with payoff $T$ between a minimizing Searcher and a maximizing Hider.
Mixed strategies are required by both players. This article adopts an
alternative search paradigm recently introduced by the authors in Alpern and
Lidbetter (2013) as ``expanding search'', because the regions $S\left(
t\right) $ that have been searched by time $t$ ``expand'' continuously over
time. This models the actual searches that are often carried out to search
for missing persons, escaped convicts, terrorists or lost airplanes. 

More formally, the sets $S\left( t\right)$, which form a pure strategy for
the searcher, start with $S\left( 0\right) =O$ and expand at unit rate of
total arc length $\lambda$, so that $\lambda \left( S\left( t\right)
\right) =t$. The search plan finishes at time $\mu =\lambda \left( Q\right) $
(thus $\mu $ is the total length of the network) because no backtracking is
required, unlike the case for pathwise search. The actual search of course
finishes at the capture time $T$, when the search region 
$S\left( T\right)$ first contains $H$. As with pathwise search, this produces a zero-sum game
with payoff $T$, which we call the {\em expanding search game}. We
concentrate on expanding searches that have a simpler structure, namely,
they consist of a sequence of arcs $a_{i},~i=1,\ldots ,k$, which cover $Q$
and where the back of each arc  $a_{i}$ lies in one of the earlier arcs.
(The direction of each arc is not prescribed in $Q$, which is an undirected
network, but simply denotes the chosen direction of search.) In our earlier
paper, it was explained how the notion of expanding search is also
applicable to non-search problems such as mining in which the cost of moving
excavation equipment through already mined trails is negligible compared to
the cost of moving the equipment by excavating new areas. Minimizing the
cost of reaching (finding) a randomly chosen piece of unmined coal is
equivalent to minimizing the mean time that kilograms of coal are mined and
ready for sale. 

\subsection{Main Results}

The expanding search game on a network was introduced in \cite{AL13}, but that paper was primarily concerned with the Bayesian
problem of minimizing the expected capture time using expanding search
against a {\em known} hider distribution. The game theoretic version only
considered two special types of networks: trees and 2-arc-connected
networks. Here we consider general networks. 

An important graph theoretic tool that we will use is the so called {\em
	bridge-block decomposition}. The arcs of a network $Q$ can be partitioned
into {\em bridges} (arcs whose removal disconnects the network) and 
{\em blocks} (the connected components of what remains after the bridges
are removed). When the blocks are each shrunk to a point, the bridges form a
tree $Q^t$. We use the parameter $r$ to denote the fraction of the total length of $Q$ that is taken up by bridges, so that when $r=0$, the network is $2$-arc-connected and when $r=1$ the network is a tree. We call $r$ the {\em bridge ratio}.

We present two new classes of mixed search strategies. In
Section \ref{sec:approx} we present the {\em block-optimal search strategy}, denoted $%
\beta$, which is optimal when there are no bridges. We show in Theorem~\ref{thm:approx} that for any network, $T( \beta ) \leq (
( 1+\sqrt{2}) /2) V(Q)$, where the value $%
V\left( Q\right) $ is the minimax search time. Note that $( 1+\sqrt{2}) /2\simeq 1.2$. We also show that
for a network with bridge ratio $r$, the expected search time $T\left( \beta \right) $ of
the block-optimal strategy satisfies the inequality $T(\beta) \le (1+r)V(Q)$, so that it performs well for ``block-like'' networks. 

In Section~\ref{sec:treelike} we present
the {\em bridge-optimal search strategy}, denoted $\gamma$, which is
optimal when there are no blocks (so $Q$ is a tree), and is based on
depth-first pure searches. Theorem~\ref{thm:gamma} shows that for any network $Q$ we have $%
T\left( \gamma \right) \leq \left( 2/\left( 1+r^{2}\right) \right) V\left(
Q\right)$. This estimate is useful when $r$ is
close to 1, so that $Q$ is ``tree-like''.

Section~\ref{sec:defs} gives a formal definition of the expanding search game. Section~\ref{sec:Bayesian}
presents results needed later on the optimal expanding search against a 
{\em known} hider distribution and begins an analysis of the  {\em
	circle-with-spike network}. Section~\ref{sec:exp-game} reviews earlier results from \cite{AL13} on the expanding search game when the network is a tree
or 2-arc-connected. Section~\ref{subsec:circle-spike-game} completes the analysis of the {\em
	circle-with-spike network }which we began in Section~\ref{sec:Bayesian}. Sections~\ref{sec:approx} and~\ref{sec:treelike}
give our main results which we have already discussed. Section~\ref{sec:conclusion} contains our
conclusions and suggestions for future work. 

\subsection{Related Literature on Network Search Games}

The use of pathwise search to minimize the time to find a hidden object on a
network was first proposed by \cite{Isaacs}. Subsequently \cite{Gal79}
analyzed such games on general networks and gave a complete solution for
trees. For trees and Eulerian networks the so-called Random Chinese Postman
Tour (RCPT), consisting of an equiprobable mixture of a minimal length
(Chinese Postman) tour and its time-reversed tour, is an optimal mixed
strategy. \cite{RP93} identified a larger class (including
trees) of networks where the RCPT is optimal. \cite{Gal00} then showed that
the largest class of networks where RCPT is optimal are the weakly Eulerian
networks. The difficult (and non-weakly Eulerian) network consisting of two
nodes (one the start node) connected by an odd number of equal length arcs
was solved by \cite{Pavlovic}. 

In addition to generalizing the classes of networks which could be solved,
other variations on the basic model have been proposed. An algorithmic
approach to the problem was given by \cite{Anderson}. Search
games on {\em windy} networks, where the times to traverse arcs from
either direction need not be the same, were introduced by \cite{Alpern2010} and
further studied by \cite{AL14}. The assumption that the game
ends not when the hider is found, but when he is brought back to the start
node, is studied in \cite{Alpern11a}. The case of expanding search with multiple hiders was solved
in some cases by \cite{Lid13a}. Search games on lattices were studied by \cite{Zoroa}. The possibility that the Searcher might have to choose between two or more speeds of search was considered by \cite{AL15}.
The case where simply reaching the hidden object is not enough to find it is
considered by \cite{BK15}, who posit a cost for searching,
in addition to traveling. Surveys of search games on networks can be found
in \cite{Garnaev}, \cite{Gal05}, \cite{Alpern11b}, \cite{Lid13b} and \cite{Hohzaki}. Network search is also related to patrolling a network to guard against an attack, as in \cite{AMP} and \cite{Lin}, for example.

The expanding search paradigm, as introduced by the authors in \cite{AL13}, has already received considerable attention in the literature. \cite{Shechter} adopt the expanding search paradigm in the constrained version of their discrete search model. \cite{Fokkink} generalize the concept of expanding search for a Hider located on one of the nodes of a tree by introducing a search model with a submodular cost function. Expanding search has provoked interest in several other areas: the search for beacons from lost airplanes in \cite{Eckman}, evolutionary game theory in \cite{Kolokoltsov} and \cite{Liu}, contract scheduling in \cite{Angelopoulos} and predator search for prey in \cite{Morice}. 

\section{Formal Definition of Expanding Search}

\label{sec:defs}

We start by repeating the formal definition of expanding search of a network given in
\cite{AL13}. Suppose $Q$ is a network with distinguished starting node $O$, which we call the {\em root}. Each arc
$a$ has a given length $\lambda\left(  a\right)$, and moreover we write
$\lambda(A)$ for the measure (total length) of any subset $A$ of $Q$, with the
total measure of $Q$ denoted  $\mu=\lambda\left(  Q\right)  $.  An expanding
search of $Q,O$ is a nested family of connected subsets of $Q$ that starts
from $O$ and increases in measure at unit speed until filling the whole of
$Q$. More formally, we use the definition from
\cite{AL13}:

\begin{definition}
	An \emph{expanding search} $S$ on a network $Q$ with root $O$ is a nested family of connected closed sets $S(t)$ for $0 \le t \le \mu$, which satisfy
	\begin{enumerate}
		\item[(i)] $S(0)=\{O\}$ and $S(\mu)=Q$,
		\item[(ii)] $S(t')\subset S(t) $ for $t'<t$, and
		\item[(iii)] $\lambda(S(t))=t$ for all $t$.
	\end{enumerate}
	We denote the set of all expanding searches of $Q$ by $\mathcal{S}=\mathcal{S}(Q)$.
\end{definition}

The first condition says that the search starts at the root and is exhaustive,
the second condition says that the sets are nested and the third says that
they increase at a unit rate.

This is the most general definition of an expanding search, but in practice we
will mostly be concerned with a class of expanding searches called
\emph{pointwise expanding searches}. A pointwise expanding search can be
thought of as a sequence of unit speed paths, each one of them beginning at a
point that has already been discovered. More formally:

\begin{definition}
	A \emph{pointwise expanding search} of a rooted network $Q,O$ is a piecewise continuous function $P:[0,\mu]\to Q$ made up of a sequence $P_1,\ldots,P_k$ of
	unit speed paths $P_{i}:[t_{i-1},t_{i}]\rightarrow Q$ with $ 0=t_0<t_1<\ldots<t_k=\mu$  such that
	\begin{enumerate}
		\item[(i)] $P_{1}(t_{0})=O$,
		\item[(ii)] $\cup _{i=1}^k P_{i}([t_{i-1},t_{i}])=Q$, and
		\item[(iii)] $P_j(t_{j-1}) \in \cup _{i=1}^{j-1} P_{i}([t_{i-1},t_{i}])$ for all $j=2,\ldots,k$.
	\end{enumerate}
\end{definition}

Note that any such sequence $\{P_i\}$ induces an expanding search $S$ where $S(t)=P([0,t])$. The first condition says that the first path starts at the root, the second
condition says the paths are an exhaustive search of $Q$ and the third says
that each path starts from a point already covered by a previous path. Not all expanding searches are not pointwise expanding searches, in particular
any expanding search that moves along more than one arc at the same time. However, in
\cite{AL13} the authors showed that the set of pointwise expanding searches is
dense in $\mathcal{S}$, so that every expanding search can be approximated to
an arbitrary degree of accuracy by a pointwise search. We will therefore use
this definition of expanding search for the rest of the paper, and will use
the terms ``expanding search'' and ``pointwise search'' interchangeably.

An important concept we will use is that of \emph{search density}, which
occurs in much previous research in search games. Suppose a Hider is located
on a network $Q$ according to some fixed probability distribution $\nu$. Then
the \emph{search density}, or simply \emph{density} $\rho=\rho\left(
A\right)  $ of a subset $A \subset Q$ is given by $\rho\left(  A\right)
=\nu\left(  A\right)  /\lambda\left(  A\right) $. When disjoint regions can be
searched in either order, it is better to search the region of highest density
first, as shown by the following simple lemma. The proof of a more general version of the lemma can be found in \cite{AL14}.

\begin{lemma}
	\label{searchdensity} Suppose the Hider is located on a rooted network $Q$ according to some given distribution $\nu$ and let $A$ and $B$ be connected regions of $Q$ that meet at
	a single point $x$. Let $S_{AB}$ and $S_{BA}$ be two expanding searches of $%
	Q $ which are the same except that, on reaching $x$, the search $S_{AB}$ follows the
	sequence $AB$ while $S_{BA}$ follows $BA$. Then if $\rho(A) \ge \rho(B)$, 
	\[
	T(S_{AB},\nu) \le T(S_{BA},\nu).
	\]
    If $A$ and $B$ have the same search density, then the regions $A$ and $B$ can be searched in either order.
\end{lemma}

\section{Searching Networks with a Known Hider Distribution}

\label{sec:Bayesian}

Before discussing the expanding search {\em game}, we give some simple results on
searching for a Hider located on a network according to a {\em known} probability
distribution. These results will be useful in our later analysis.

For a given network $Q$ with root $O$, suppose the Hider is located at some
point $H$ on the network. For a given expanding search $P$ we will write
$T(P,H)$ to denote the first time that $H$ is \textquotedblleft
discovered\textquotedblright\ by $P$. That is, $T(P,H)=\min\{t\geq0:P(t)=H\}$.
This was shown to be well defined in
\cite{AL13}. We call $T(P,H)$ the \emph{search time}. If the Hider is located
on $Q$ according to some probability distribution $\nu$, we denote the
\emph{expected} search time of a search $P$ by $T(P,\nu)$ (or later, in the expanding search game, by $T\left(
p,\nu\right)  $ if the Searcher is adopting a mixed search strategy $p$).

For a given Hider distribution $\nu$ we are interested in the problem of
finding the expanding search that minimizes the expected search time. We call such a search {\em optimal}. In
\cite{AL13} the authors showed that for a rooted tree with a Hider located on
it according to a known distribution, the optimal search begins with the
rooted subtree of maximum density (assuming it exists). More formally the
theorem is:

\begin{thm}[Theorem 14 of
	\cite{AL13}]
	\label{old theorem}Let the Hider $H$ be hidden according to a known
	distribution $\nu $ on a rooted tree $Q$ and suppose there is a unique
	rooted subtree $A$ of maximum density. Then there is an optimal expanding
	search $S$ which begins by searching $A$. That is, $S\left( \lambda \left(
	A\right) \right) =A$.
\end{thm}

When the Hider distribution $\nu$ is concentrated on the nodes of $Q$, this
theorem gives a simple algorithm for computing an optimal expanding search of
$Q$. The theorem is also true when $\nu$ is a continuous distribution.

We will later use a result that says if $A$ is a component of $Q$ on which the
Hider is hidden uniformly, and $A$ is connected to $Q$ at only one point, then
there is an optimal search that searches the whole of $A$ at once. (Formally,
when we say the Hider is hidden uniformly on $A$, we mean that given he is on
$A$, the probability he is located within some measurable subset $X$ of $A$ is
proportional to the measure of $X$.)

\begin{lemma}
	\label{uniform}Suppose $Q$ is composed of two subnetworks $A$ and $B$ which
	meet at a single point. Suppose a Hider is located on $Q$ according to some
	distribution which is uniform on $A$. Then there is an optimal expanding
	search that searches $A$ without interruption.
\end{lemma}

\proof{Proof.} Any expanding search of $Q$ must search $A$ in a finite number
of (closed) disjoint time intervals. Let $P$ be an optimal expanding search of
$Q$ that searches $A$ in a minimal number $m$ of disjoint time intervals. If
$m=1$ then the lemma is true, so suppose $m \ge2$ and we will derive a contradiction.

Let $A_{1}$ and $A_{2}$ be the subsets of $A$ searched in the first two
disjoint time intervals that $A$ is searched, and let $C$ be the subset of $B$
that is searched between $A_{1}$ and $A_{2}$. Note that since $A$ and $B$ meet
at a single point, it is possible to search $C$ at any point in time after $P$
starts searching $A$. By Lemma \ref{searchdensity}, the search density of $C$
is no more than that of $A_{1}$, otherwise the expected search time could be
reduced by searching $C$ before $A_{1}$. By a similar argument, the search
density of $A_{2}$ is no more than that of $C$. In other words,
\[
\rho(A_{1}) \ge\rho(C) \ge\rho(A_{2}).
\]

But since the Hider is hidden uniformly on $A$, the search density of $A_{1}$
and $A_{2}$ are the same so both of the inequalities above hold with equality.
Hence, by Lemma \ref{searchdensity}, the search $P^{\prime}$ that is the same
as $P$ except that $A_{2}$ and $C$ are searched the other way around is also
optimal. But $P^{\prime}$ searches $A$ in $m-1$ disjoint time intervals, a
contradiction. 
\endproof

Every expanding search $P$ leads to a \emph{search tree} $Q_{P}$
obtained from $Q$ by cutting it at certain points $x\in Q$ whose
removal leaves a tree. Roughly speaking, these points $x$ are those which are
reached from more than one direction by the search $P$. If a search $P$ is
optimal, it has to be optimal for the corresponding distribution on the tree
$Q_{P},$ so $P$ can be found by applying Theorem~\ref{old theorem} to all such trees.

\begin{figure}[th]
	\begin{center}
		\includegraphics[scale=0.5]{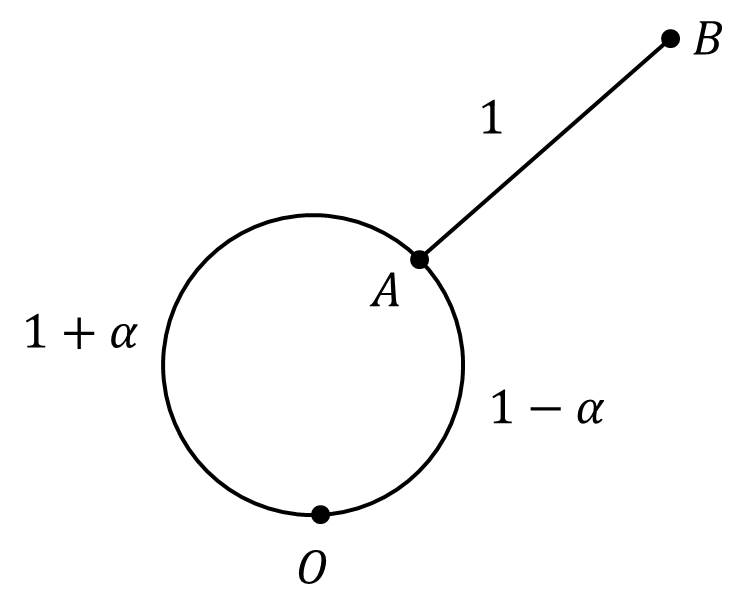}
	\end{center}
	\caption{The circle-with-spike network $Q^{CS}$.}%
	\label{fig:circle-spike}%
\end{figure}

Consider the circle-with-spike network $Q^{CS}$ drawn in Figure
\ref{fig:circle-spike}. It consists of a circle of length $2$, to which a unit
length line segment is attached at a clockwise distance $1+\alpha$ from the
root $O,$ $0\leq\alpha<1$. Consider Hider distributions $\nu_{p}$ having atoms
of weight $\nu_{p}\left(  B\right)  =1-p$ at leaf node $B$ and a uniform
distribution of total weight $p$ on the clockwise arc from $O$ to $A$. For
fixed $\alpha$ and $p$, we determine the optimal expanding search. Every
expanding search produces a subtree of $Q^{CS},$ so our method is first to
determine these subtrees and then to find the optimal search corresponding to
that subtree using Theorem \ref{old theorem}.

Since the Hider distribution on the two arcs of the circle is uniform, by
Lemma \ref{uniform}, we need only consider making cuts at the nodes.
Hence there are only two trees we need consider, determined by the initial arc
of the search. For the search beginning by going clockwise, $P_{+}$, the tree
obtained by disconnecting the counterclockwise arc at $O$ is optimally
searched by going to $B$ after reaching $A$, as the Hider is certain to be
there if not already found. So we have
\[
T\left(  P_{+},\nu_{p}\right)  =p\left(  \frac{1+\alpha}{2}\right)  +\left(
1-p\right)  \left(  2+\alpha\right)  .
\]
For the search that begins by going counterclockwise, $P_{-},$ the tree is
obtained by disconnecting the clockwise arc from $O$ (or equivalently,
disconnecting it from $A$). If this arc has higher density than the one going
to $B$, it would have been better to use $P_{+}$, so we may assume that going
to $B$ has higher density, and that
\[
T\left(  P_{-},\nu_{p}\right)  =p\left(  2-\alpha+\frac{1+\alpha}{2}\right)
+\left(  1-p\right)  \left(  2-\alpha\right)  .
\]
Equating these times, we see that the search $P_{+}$ is better for $p<2\alpha/\left(
\alpha+2\right)$, the search $P_{-}$ is better for $p>2\alpha/\left(  \alpha+2\right)
$ and they give the same times for $\bar{p}=2\alpha/\left(  \alpha+2\right)
.$ So in particular the least search time against $\nu_{\bar{p}}$ is
\[
\min_{P} T\left(P,  \nu_{\bar{p}}\right)  =T\left(  P_{+},\nu_{\bar{p}}\right)
=\frac{\alpha+4}{\alpha+2}.
\]
This analysis proves only that $\nu_{\bar{p}}$ is the hardest Hider
distribution to find within the family $\nu_{p},$ but later in Section
\ref{subsec:circle-spike-game} when we consider a game theoretic analysis, we
will establish that in fact it is the hardest Hider distribution to search,
without any restrictions. We note that if $\alpha\leq2/3$ then putting the
uniform distribution over the long arc onto its center does not reduce its
minimum search time. However if this is done when $\alpha> 2/3$, the Searcher
can replace the strategy that begins on the long arc by one which goes to its
center and then traces out the short arc from the root, saving some time as
this gets to $A$ faster.

\section{Previous Results on the Expanding Search Game}
\label{sec:exp-game}
We now assume that the Hider distribution $\nu$ is not known to the Searcher.
In this case we consider the problem of finding the mixed Searcher strategy
(probability measure over expanding searches) which minimizes the expected
search time in the worst case. An equivalent problem, which we prefer to
adopt, is the zero-sum Expanding Search Game $\Gamma\left(  Q\right)  $. Here
the maximizing Hider picks a location $H$ in $Q$, the minimizing Searcher
picks an expanding search $P$ and the payoff is the search
(capture) time $T\left(  P,H\right)  $. The analogous game, $\Gamma^{p}(Q)$
where the Searcher picks a unit speed path on $Q$ has been well studied, and
we call this model of search \emph{pathwise search}, denoting the value of the
analogous \emph{pathwise search game} by $V^{p}(Q)$. We note that $V\leq
V^{p}$.

In
\cite{AL13}, the authors showed that for any network $Q$ the expanding search
game has a value, $V=V(Q)$, and they solved the game in the cases that $Q$ is
a tree and that $Q$ is $2$-arc-connected (that is $Q$ cannot be disconnected
by the removal of fewer than $2$ arcs). In general the Searcher is not
restricted to using a pointwise search, but in fact in the solutions for these
classes of networks the Searcher always randomizes between pointwise searches.

We first present the solution of the game for trees. In this case, the Hider's optimal distribution is
concentrated on the leaf nodes of the network, since all others points are
dominated. We will only consider binary trees (that is trees with maximum degree at most
3), since any tree can be transformed into a binary tree by adding arcs of
arbitrarily small length. For a branch node $x$ of a rooted tree, $Q,O$ with
outward arcs $a$ and $b$, we denote the branches starting with $a$ and $b$ by
$Q_{a}$ and $Q_{b}$, respectively, and their union by $Q_{x}$. We also denote
the length of $Q_{a}$ and $Q_{b}$ by $\mu_{a}$ and $\mu
_{b}$, and write $\mu_{x}=\mu_{a}+\mu_{b}$ for the length of $Q_{x} = Q_{a} \cup Q_{b}$.

\begin{definition}
	\label{EBD}
	Let $Q,O$ be a rooted, variable speed tree. Let the \emph{Equal Branch Density} (\emph{EBD}) distribution, $e$ be the unique probability distribution on the leaf nodes $\mathcal{L}(Q)$ of $Q$ such that at any branch node $x$ of $Q$, all the branches rooted at $x$ have the same search density.
\end{definition}

For a branch node $x$ we denote the EBD distribution on the subnetwork $Q_{x}$
by $e_{x}$ and similarly, for an outward arc $a$ of $x$ we denote the EBD
distribution on $Q_{a}$ by $e_{a}$.

In order to describe the optimal Searcher strategy we need to define a quantity $D(Q)$, which is the average distance from the root of $Q$ to its leaf nodes, weighted with respect to the EBD distribution.

\begin{definition}
	For a rooted tree, $Q,O$, the quantity $D=D(Q)$ is defined by
	\[
	D(Q) = \sum_{i \in \mathcal{L}(Q)}e(i)d(O,i),
	\]
	where $d(O,i)$ is the length of the path from $O$ to $i$ in $Q$.
\end{definition}

\begin{definition}
	We define the \emph{biased depth-first strategy} for the Searcher as the mix of depth-first searches such that at any given branch node $x$ with branches $Q_a$ and $Q_b$, on encountering $x$ for the first time, the Searcher searches the whole of $Q_a$ first with probability $p(a)$ given by
	\[
	p(a) = \frac{1}{2} + \frac{1}{2\mu_x}(D(Q_a)-D(Q_b)).
	\]
\end{definition}

In
\cite{AL14}, the authors prove that the solution of the expanding search game on a tree is as follows.

\begin{thm}[Theorems 18 and 19 of
	\cite{AL13}]
	\label{treegame}
	The value $V$ of the expanding search game played on a rooted tree $Q,O$ is given by
	\[
	V=\frac{1}{2}(\mu+D).
	\]
	The EBD distribution is optimal for the Hider and the biased depth-first strategy is optimal for the Searcher.
\end{thm}

We also state the solution of the expanding search game for 2-arc-connected
networks. We first define a \emph{reversible expanding search} as an expanding
search whose time-reverse is also an expanding search. Not all expanding
searches are reversible, but it was shown in
\cite{AL13} that a network admits a reversible expanding search if and
only if it is 2-arc-connected. This gives rise to the following theorem from
the same paper.

\begin{thm}[Theorem 20 and Corollary 22 of
	\cite{AL13}]
	\label{thm:2arcconnected}
	A rooted tree $Q,O$ is 2-arc-connected if and only if it admits a reversible expanding search. If $Q$ is 2-arc-connected, then the value $V$ of the expanding search game is $\mu/2$. An equiprobable choice of some reversible expanding search and its reverse is optimal for the Searcher; the uniform distribution is optimal for the Hider.
\end{thm}

\section{Search Game on Circle-with-Spike Network}

\label{subsec:circle-spike-game}

We continue our analysis of expanding search games by considering the
circle-with-spike networks $Q^{CS}$ of Section \ref{sec:Bayesian}. The
analysis given there of the Hider distributions $v_{\bar{p}}$ shows that
$V=V\left(  Q^{CS}\right)  \geq\left(  \alpha+4\right)  /\left(
\alpha+2\right) $. We now give a mixed Searcher strategy which shows that
equality holds.

The Searcher can find the Hider in time at most $\left(  4+\alpha\right)
/\left(  2+\alpha\right)  $ by using the mixed strategy $\sigma_{\alpha}$ of
choosing searches $(P_{1},P_{2},P_{3})$ with probabilities $\left(  \frac
{1}{2},\frac{1}{2(2+\alpha)},\frac{1+\alpha}{2(2+\alpha)}\right)  $, where the
$P_{i}$ are given as follows: $P_{1}$ travels anticlockwise along the circle
to $A$, then goes to $B$, and finally traverses the remaining arc from $A$ to
$O$; $P_{2}$ is the same as $S_{1}$, but traverses the last arc from $O$ to
$A$; $P_{3}$ travels clockwise along around the circle to $A$, then goes to
$B$, and finally traverses the remaining arc from $A$ to $O$. Note that hiding
anywhere in the arc $AB$ is dominated by hiding at $B$ and that against
$\sigma_{\alpha}$ hiding on the short circular arc is dominated by $B.$

If the Hider is a clockwise distance $x<1+\alpha$ from $O$,\ the expected
search time is%
\[
\frac{1}{2}\left(  (2-\alpha)+(1+\alpha-x)\right)  +\frac{1}{2(2+\alpha
	)}(2-\alpha+x)+\frac{1+\alpha}{2(2+\alpha)}\left(  x\right)  =\frac{4+\alpha
}{2+\alpha}\text{.}
\]
If the Hider at at $B,$ the expected search time is%
\[
\left(  \frac{1}{2}+\frac{1}{2(2+\alpha)}\right)  \left(  2-\alpha\right)
+\left(  \frac{1+\alpha}{2(2+\alpha)}\right)  \left(  2+\alpha\right)
=\frac{4+\alpha}{2+\alpha}.
\]
So we have shown the following. \ 

\begin{thm}
	For the circle-with-spike networks $Q^{CS}$ the value of the expanding
	search game is given by $V=\left( 4+\alpha \right) /\left( 2+\alpha \right)
	. $ The optimal Hider strategy is $\nu _{\bar{p}},$ $\bar{p}=2\alpha /\left(
	\alpha +2\right) ;$ the optimal Searcher strategy is $\sigma _{\alpha }.$
\end{thm}

It is useful to compare this result with the solution of the pathwise search
game on the same network $Q^{CS}$. These networks are weakly Eulerian. Roughly
this means they consist of disjoint Eulerian networks which when shrunk to a
point leave a tree. (In particular here such a shrinkage leads to the tree
$AB,$ with the circle shrunk to root $A$.) The work of
\cite{RP93} - because these networks are also weakly cyclic - and
\cite{Gal00} show that the pathwise search value is $V^{p}=\bar{\mu}/2=2$, where
$\bar{\mu}$ denotes the shortest tour time, which here is the length of the
circle plus twice the length of the spike. This is the same as the expanding
search value when $\alpha=0,$ that is, when the point $A$ is antidpodal to the
root. Otherwise, the expanding search value is strictly smaller. In pathwise
search, the Hider optimally hides uniformly over the circle regardless of the
location of the root on the circle.

For weakly Eulerian networks, the pathwise search value $V^{p}$ is independent
of the location of the root. This is no longer true for such networks in
expanding search. In this example, where the location of the root depends on
the parameter $\alpha$, we have shown that the value $V_{\alpha}=\left(
4+\alpha\right)  /\left(  2+\alpha\right)  $ is also dependent on $\alpha$.

\section{The Block-Optimal Search Strategy}

\label{sec:approx}

In
\cite{AL13}, optimal expanding search strategies were found for networks which
are 2-arc-connected or trees. In the present paper optimal expanding search
strategies have also been found for the circle-with-spike networks. Since optimal
strategies for searching other networks are not known, it is of use to find a
class of strategies which can be calculated for any rooted network, and are
``approximately optimal''. We define a version of approximate optimality for the Searcher
which is multiplicative in nature. We say that a family $s_{Q}$ of mixed expanding search strategies for rooted
networks $Q$ is {\em $\alpha$-factor optimal} if for every $Q$, we have%
\[
T(s_{Q},H) \leq\alpha V(Q)\text{, for all }H\in Q,
\]
where $\alpha\ge1$ is a constant independent of $Q$. Writing $T(s_{Q})$ for
the maximum value $T(s_{Q},H)$ takes over all $H$, we can equivalently write
that the family $s_{Q}$ is $\alpha$-factor optimal if $T(s_{Q}) \le\alpha
V(Q)$. The closer $\alpha$ is to $1$, the better the approximation given by
the search family $s_{Q}$. In the language of approximation algorithms, the
notion of an $\alpha$-factor optimal strategy is akin to the notion of an
$\alpha$-approximate algorithm for computing the value of the game.

In this section we present such a family of search strategies, called
\emph{bridge-optimal strategies} $\beta=\beta_{Q}$, which are $1.2$%
-factor optimal. The name comes from the fact that these strategies use the
well-known ``bridge-block'' decomposition of an arbitrary network, as defined in Subsection~\ref{sec:bridge-block}. The
multiplicative constant $\alpha$ is actually $\left(  1+\sqrt{2}\right)  /2$
which is approximately equal to $1.207$.

It turns out that the range of factor constants $\alpha$ that we need to
consider goes from $1$ to $2$. To see that a constant of $\alpha=2$ is of no
value consider any class $s_{Q}$ of expanding search strategies. The
definition of an expanding search shows that for any point $H\in Q$ we have
$T(s_{Q},H) \leq\mu$, the total length of $Q$. Since the uniform hiding
strategy ensures an expected search time of at least $\mu/2$, we know that
$V(Q) \geq\mu/2$. Consequently, for any family $s_Q$, we have
\[
T(s_{Q},H) \leq\mu\leq2(\mu/2) \leq2V(Q) \text{, for all }H\in Q,
\]
and thus \emph{any} expanding search family $s_{Q}$ is $2$-factor optimal.

\subsection{The Bridge-Block Decomposition}
\label{sec:bridge-block}

Let $Q$ be a connected network. An arc of $Q$ is called a \emph{bridge} if
removing it (but not its end nodes) disconnects the graph, or equivalently,
if it is not contained in any cycle. The components of $Q$
after removing its bridges are called the \emph{blocks}. Note that the blocks are $2$-arc-connected. We denote the set of
bridge-points (points in bridges) of $Q$ by $Q_{1}$ and the set of
block-points (points in blocks) by $Q_{2}$. The connected network $Q^t$
obtained from $Q$ by shrinking each block to a point (node) is a tree whose
arcs are identical to the arcs of $Q_1$. We call $Q^t$ the \emph{bridge
	tree} of $Q$. The nodes of $Q^t$ are of two types: the \emph{new nodes}
correspond to the blocks of $Q$ and the \emph{original nodes} correspond to
nodes of $Q$ which are incident only to bridges in $Q$.

\begin{figure}[th]
	\begin{center}
		\includegraphics[scale=0.5]{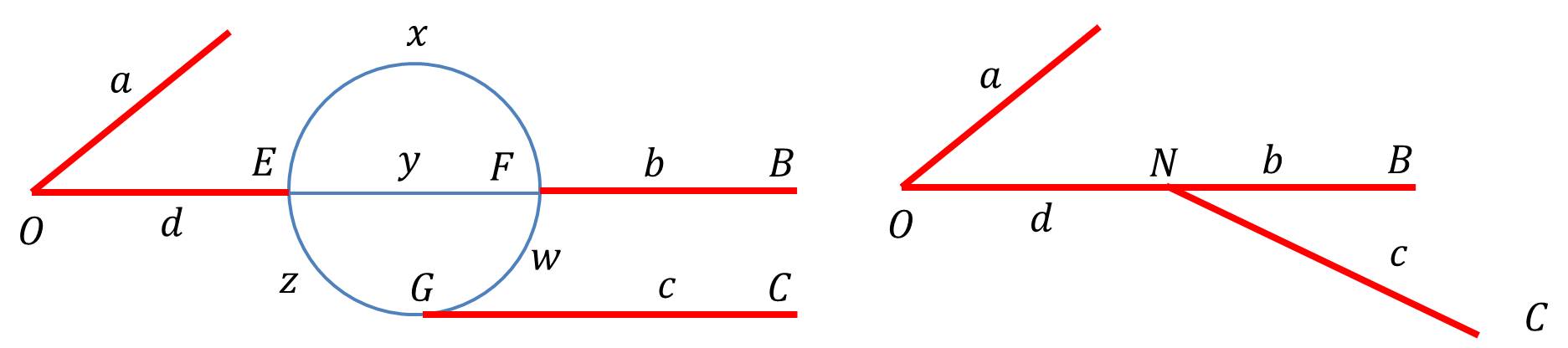}
	\end{center}
	\caption{A network $\bar{Q}$ (left) and its bridge tree $\bar{Q}^{t
		}$.}%
	\label{fig:bridge-block}%
\end{figure}

These concepts are illustrated in Figure \ref{fig:bridge-block}. The network
$\bar{Q}$ on the left has four bridges: $a,b,c,d$. It has a single block made
up of arcs $x,y,z,w$. Its bridge tree $\bar{Q}^t$ consists of the 4 arcs
$a,b,c,d$. It has a new node $N$ which corresponds to the block of $\bar{Q}$.

 For any point $H$ in
 $Q$, its height $\pi(H)$ is
 defined as the distance from the root to its corresponding point in $Q^t%
 $, so that all points in a given block have the same height. The notion of height is illustrated in Table~\ref{tab:Qbar} and Figure~\ref{fig:bridge-block tree}, where we assign lengths to the arcs of $\bar{Q}$.

We also define the {\em bridge ratio} $r=r(Q)$ of a network $Q$ as the fraction of the total length of $Q$ that consists of bridges. So, writing $\mu_1$ for the total length of all the bridges in $Q$, the bridge ratio is given by $r=\mu_1/\mu$.

\subsection{Definition of the Block-Optimal Search Strategy $\beta$}

We now define the block-optimal search strategy $\beta$ as an equiprobable mixture of two
expanding searches we call $S_{1}$ and $S_{2}$. The second one, $S_{2}$, will
be a sort of reverse path of $S_{1}$ so we concentrate on describing $S_{1}$.
The specification of $S_{1}$ will not depend on arc lengths, so we leave that
until later. To specify $S_{1}$ we fix on each of its blocks a particular
reversible expanding search (as defined in Section~\ref{sec:exp-game}). For example, on the pictured network
$\bar{Q}$ we orient each arc to the right and denote a traversal to the left
with a prime. So one reversible expanding search of its block is given by
$x,y^{\prime},w^{\prime},z^{\prime}$ with reverse path $z,w,y,x^{\prime}$
(reverse order and direction). The fact that every block has a reversible
expanding search follows from Theorem \ref{thm:2arcconnected}. The expanding
search $S_{1}$ follows the first of these two reversible expanding searches ($x,y^{\prime},w^{\prime
},z^{\prime}) $ on each block. However when reaching a cut-node such as $F$ or
$G$ which leads out of the block, it exhaustively searches that component of
the network before returning to node on which it left. When reaching a node of
$Q $ which is not in a block (in $\bar{Q}$ the only such node is the start
node $O$) $S_{1}$ can leave via any arc; the reverse ordering will be chosen
by $S_{2}.$

We describe the construction of $S_{1}$ for the pictured network $\bar{Q}.$ We
begin by fixing a reversible expanding search at its only block:
($x,y^{\prime},w^{\prime},z^{\prime})$. Now we start at $O$ and can choose
either of the arcs $a$ or $d.$ We choose $d$ (this will mean that $S_{2}$
begins with $a$). We arrive at node $E$ which is in the block so we travel
according to our fixed reversible path until we reach a node which leaves the
block. So we must continue with $d,x$ which arrives at the node $F.$ Here we
leave the block by choosing $d,x,b.$ Since the component of $Q-F$ is now
exhaustively searched we continue searching the block from $F,$ continuing
$d,x,b,y^{\prime},w^{\prime}$ which brings us to another node, $G,$ which
leads out of the block. Here we continue $d,x,b,y^{\prime}w^{\prime
},c,z^{\prime}.$ Now we have exhaustively searched the part of $\bar{Q}$
stemming from the arc $d.$ We now follow the same procedure on the part of
$\bar{Q}$ stemming from $a,$ which in this case is just $a,$ leading to the
construction
\[
S_{1}=d,x,b,y^{\prime},w^{\prime},c,z^{\prime},a.
\]
To construct $S_{2}$ we follow the reverse expanding search on each block,
here $z,w,y,x^{\prime}.$ We also take the opposite ordering on the non-new
branch nodes of the bridge tree, that is, we start with $a$ rather than $d.$
So we start with $a,d.$ Now we follow the reverse expanding search on the
block, beginning with $a,d,z.$ Then we search the entire portion of the
network that begins with arc $c,$ which in this case is just $c$ itself,
obtaining $a,d,z,c.$ We continue searching the block with $a,d,z,c,w.$ Since
$b$ came before $y$ (actually $y^{\prime}$) in $S_{1}$, we continue in the
opposite order in $S_{2},$ with $a,d,z,c,w,y,b$ and then finish with
$x^{\prime},$ obtaining the expanding search
\[
S_{2}=a,d,z,c,w,y,b,x^{\prime}.
\]
We now explore the performance of $\beta=\left(  1/2\right)  S_{1}+\left(
1/2\right)  S_{2}$ against undominated hider strategies. A hider strategy $H$
is undominated if it lies on an arc of a block or on one of the leaf nodes in
the bridge tree. Node that arcs in blocks are traversed in opposite directions
by $S_{1}$ and $S_{2},$ so it does not matter where the Hider is situated in
such an arc. We take the lengths of $a,d,x,y,b$ as $2$; $c$ as 3; $z$ and $w$
as $1$. So $\mu=15$. Table \ref{tab:Qbar} gives the times take for $S_{1}$ and
$S_{2}$ to reach a Hider placed at one of the leaf nodes $A,B,C$ or at the
center of $x,y,w,z$. \begin{table}[tbh]
	\caption{Search times for the network $\bar{Q}$.}%
	\label{tab:Qbar}
	\begin{center}%
		\begin{tabular}
			[c]{|l|l|l|l|l|l|}\hline
			& $T(S_{1},H)$ & $T(S_{2},H)$ & sum & height, $\pi(H)$ & $\mu+ \pi(H) = 15 +
			\pi(H)$\\\hline
			$A$ & 15 & 2 & 17 & 2 & 17\\\hline
			$B$ & 6 & 13 & 19 & 4 & 19\\\hline
			$C$ & 12 & 8 & 20 & 5 & 20\\\hline
			$x$ & 2.5 & 14.5 & 17 & 2 & 17\\\hline
			$y$ & 7 & 10 & 17 & 2 & 17\\\hline
			$z$ & 12.5 & 4.5 & 17 & 2 & 17\\\hline
			$w$ & 8.5 & 8.5 & 17 & 2 & 17\\\hline
		\end{tabular}
	\end{center}
\end{table}

\begin{figure}[th]
	\begin{center}
		\includegraphics[scale=0.5]{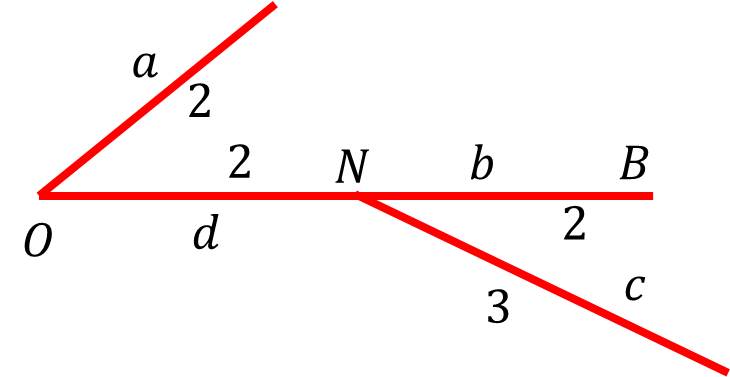}
	\end{center}
	\caption{The bridge tree $\bar{Q}^t$ with arc lengths.}%
	\label{fig:bridge-block tree}%
\end{figure}

There is an easier way to calculate the expected time (the column labeled
\textquotedblleft sum\textquotedblright\ divided by 2) for the block-optimal strategy
$\beta$ to reach a point in $\bar{Q}$ corresponding to any node of the bridge
tree $\bar{Q}^t.$ To understand this method, consider the bridge tree as
drawn in Figure~\ref{fig:bridge-block tree} with arc lengths of bridges
indicated. Note that the only hiding locations that are best responses to
$\beta$ correspond to leaf nodes or new (block) nodes of $Q^t$ (for
$\bar{Q}^t$ all nodes other than $O$ are of one of these types). Suppose
for example, that the Hider chooses $H=C.$ Then certainly the whole lengths of
arcs $d$ and $c$ must be covered before reaching $H$, so the expected search
time is at least $5$. Note that every other point will be reached before $H$
by exactly one of the pure searches $S_{1}$ and $S_{2}$. Since these
\textquotedblleft other points\textquotedblright\ have total length
$\mu-5=15-5=10$, the expected search time is $5+\left(  1/2\right)  10=10$, as
given in our table as $20/2=sum/2.$ This argument is quite general for any
network $Q$. Suppose the Hider hides at a leaf node or on a non-bridge arc
$a$, breaking it into arcs $a^{\prime}$ and $a^{\prime\prime}$ when $H$ is
considered a node. Then as in the above argument, the bridge arcs which
connect $A$ to $O$ in $Q^t$ have total length $\pi\left(  H\right)  $
and must definitely be traversed by $\beta$ before it finds $H.$ All other
arcs are traversed either once by $S_{1}$ or once by $S_{2}$ (so on average
1/2 by $\beta)$ so the expected time for $\beta$ to $H$ is given by
\begin{align*}
T(\beta,H) &  =\pi(H)+(1/2)(\mu-\pi(H))\\
&  =\frac{\mu+\pi(H)}{2},\text{ for }H\in Q_{2}\text{ or }H\text{ a leaf
	node.}%
\end{align*}
Note that for other hider places $H$ the search time may be smaller, so we
have
\[
T(\beta,H)\leq\frac{\mu+\pi(H)}{2},\text{ for all }H\in Q.
\]

Consequently taking $\pi=\pi(Q)=\max_{H\in Q}\pi(H)$ to be the \emph{height}
of $Q$, we have%
\begin{equation}
T(\beta)=\frac{\mu+\pi}{2}.\label{eq:max-time}%
\end{equation}

Note that there must be a leaf node $H$ of $Q_{1}$ such that $\pi(H)=\pi$.

\subsection{Performance of $\beta$}

We now show that the block-optimal strategy $\beta$ is $1.2$-factor optimal. To do this, we first find a lower bound on $V(Q)$. This is turn is accomplished by finding the value of the game on certain networks $Q'$ and $Q''$ which are easier to search than $Q$. 

\begin{figure}[th]
	\begin{center}
		\includegraphics[scale=0.5]{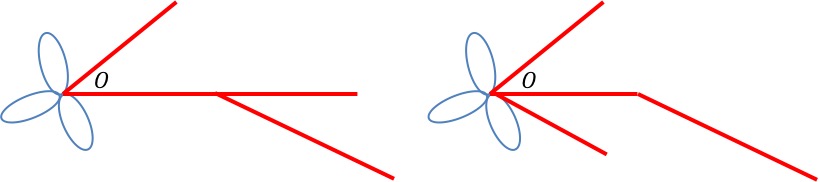}
	\end{center}
	\caption{The networks $Q'$ (left) and $Q''$ (right) depicted for the network $Q=\bar{Q}$ of Figure~\ref{fig:bridge-block}.}
	\label{fig:pruning}
\end{figure}

\begin{lemma}
	\label{lem:pruning}
	For any rooted network $Q$, the value $V(Q)$ of the game satisfies
	\[
	V(Q) \ge \frac{1}{2 \mu}(\mu^2 + \pi^2).
	\]
\end{lemma}
\proof{Proof.} 
We start by defining a new network $Q^{\prime}$ whose value is
no greater than the value of $Q$. The network $Q^{\prime}$ is obtained by
first identifying all the nodes in each block of $Q$, so that each block
is now a set of loops that meet at a single node. We then remove those loops
and reattach them at the root $O$. So $Q^{\prime}$ consists of the bridge tree $Q^t$, with some extra loops meeting at $O$. This network $Q'$ is depicted on the left of Figure~\ref{fig:pruning} for the network $Q=\bar{Q}$ of Figure~\ref{fig:bridge-block}. The value
$V(Q^{\prime})$ is lower than $V(Q)$ since any Searcher strategy that can be
used on $Q$ can also be used on $Q^{\prime}$.

We then define another new network $Q^{\prime\prime}$ whose value is, again,
no greater than the value of $Q^{\prime}$. Let $H^{*}$ be a leaf node of
$Q_{1}$ at distance $\pi$ from $O$. The network $Q^{\prime\prime}$ is obtained
from $Q^{\prime}$ by ``pruning'' the path $P$ from $O$ to $H^{*}$: all
subtrees rooted on this path are removed and reattached at $O$, so that in
$Q^{\prime\prime}$ all the vertices on $P$ except $O$ and the leaf nodes have degree
$2$. The network $Q''$ is depicted on the right of Figure~\ref{fig:pruning}. It is clear that the value $V(Q^{\prime\prime})$ is no greater than
$V(Q^{\prime})$, since any Searcher strategy on $Q^{\prime}$ can also be
executed on $Q^{\prime\prime}$.

We now derive a lower bound on $V(Q^{\prime\prime})$ (and therefore on $V(Q)$)
by giving an explicit Hider mixed strategy $h$ on $Q^{\prime\prime}$ and
calculating its minimum expected search time against any Searcher strategy.
The Hider strategy $h$ hides at $H^{\ast}$ with probability $\pi/\mu$ and
hides uniformly on the rest of $Q''$. This means that the search density of the
path $P$ from $O$ to $H^{\ast}$ is the same as the search density of the rest
of the network. By Lemma~\ref{uniform}, the two regions $P$ and $Q-P$ must be searched one
after the other, and the fact the densities are equal means it does not matter
which way around they are searched. So let $S$ be the search that first visits
$H^{\ast}$ then searches the rest of the network (it does not matter how). The
expected search time is
\begin{align*}
T(S,h) &  =\frac{\pi}{\mu}(\pi)+\left(  1-\frac{\pi}{\mu}\right)  \left(
\pi+\frac{\mu-\pi}{2}\right)  \\
&  =\frac{1}{2\mu}\left(  \mu^{2}+\pi^{2}\right)  .
\end{align*}
\endproof
By comparing Equation~(\ref{eq:max-time}) with Lemma~\ref{lem:pruning}, we can estimate the efficiency of the strategy $\beta$.
\begin{thm}
	\label{thm:approx}
	Let $Q$ be a rooted network with height $\pi$ and total length $\mu$. The maximum expected search time $T(\beta)$ of the block-optimal strategy $\beta$ on $Q$ satisfies
	\begin{align}
	T(\beta) & \leq (1+r) V(Q). \label{eq:beta-approx-r}
	\end{align}
Moreover, $\beta$ is $(1+\sqrt{2})/2$-factor optimal for any network.
\end{thm}

\proof{Proof.} 
It follows from Equation~(\ref{eq:max-time}) and Lemma~\ref{lem:pruning} that
\begin{align}
\frac{T(\beta)}{V(Q)} &  \leq\frac{\frac{1}{2}(\mu+\pi)}{\frac{1}{2\mu}%
	(\mu^{2}+\pi^{2})}  =\frac{1+\pi/\mu}{1+(\pi/\mu)^{2}}. \label{eq:lowerbound1}
\end{align}
Note that the definition of the height $\pi$ ensures that it cannot exceed the total length of the bridges, so $\pi \le \mu_1$ and hence
\[
\frac{\pi}{\mu} \le \frac{\mu_1}{\mu} \equiv r.
\]
Thus (\ref{eq:beta-approx-r}) follows from~(\ref{eq:lowerbound1}).

To show that $\beta$ is $(1+\sqrt{2})/2$-factor optimal, it is easy to check that the right-hand side of~(\ref{eq:lowerbound1}) is maximized when $\pi/\mu=\sqrt{2}-1$, and so we get
\[
\frac{T(s)}{V(Q)}\leq\frac{\sqrt{2}}{1+(\sqrt{2}-1)^{2}}=(1+\sqrt
{2})/2.
\]
\endproof
The constant $(1+\sqrt{2})/2$ is the best possible in general. For example,
suppose $Q$ is the star network consisting of one arc of length $1$ and $n$
arcs length $\sqrt{2}/n$. Then $\mu=1+\sqrt{2}$ and $\pi=1$ so by
(\ref{eq:max-time}), we have $T(\beta)=(2+\sqrt{2})/2$. But by Theorem
\ref{treegame}, the value $V(Q)$ can be calculated to be
\[
V(Q)=\frac{1}{2}(\mu+D)=\frac{2+\sqrt{2}+1/n}{1+\sqrt{2}}\rightarrow
\frac{2+\sqrt{2}}{1+\sqrt{2}}\mbox{ as }n\rightarrow\infty.
\]
Hence the ratio between $T(\beta)$ and $V(Q)$ can be arbitrarily close to
$(1+\sqrt{2})/2$.

Using~(\ref{eq:lowerbound1}), we can say a bit more precisely how well the block-optimal strategy $\beta$ performs for specific networks. For the circle-with-spike network $Q^{CS}$, the height is $\pi=1$, so $\pi/\mu=1/3$ and by~(\ref{eq:lowerbound1}), we have $T(\beta) \le 1.2 V(Q^{CS})$. The ratio of the height to the total measure of the network $\bar{Q}$ of Figure~\ref{fig:bridge-block} is also $\pi/\mu=5/15=1/3$, so we also have $T(\beta) \le 1.2 V(\bar{Q})$.

The inequality~(\ref{eq:beta-approx-r}) makes it clear that the less $Q$ is taken up with bridges, the better $\beta$ performs. In particular, it is optimal if $Q$ contains no bridges, and is therefore $2$-arc-connected. This is also obvious from the definition of $\beta$, by Theorem~\ref{thm:2arcconnected}.

\section{The Bridge-Optimal Search Strategy}
\label{sec:treelike}

In Section~\ref{sec:approx}, we found that the block-optimal strategy $\beta$ approximates the optimal search strategy within a factor of 1.2 and performs particularly well on ``block-like'' networks. In this section we will give an alternative search strategy, $\gamma$ that performs well in the case that the network mostly consists of bridges with a small proportion of the length taken up by blocks, so that the network is ``tree-like''. We call $\gamma$ the {\em bridge-optimal strategy}.

The strategy $\gamma$ is simply an adaptation of the biased depth-first strategy, described in Section~\ref{sec:exp-game}, which is optimal for tree networks. The strategy $\gamma$ follows the biased depth-first strategy on the set of bridge-points $Q_1$, as if using the optimal strategy on the bridge tree $Q^t$ of $Q$; and whenever the Searcher encounters a block in $Q_2$ for the first time, she performs an arbitrary search of the whole of it. Since the biased depth-first strategy is optimal for $Q^t$, for any point $x$ of $Q$ that lies on a bridge, the strategy $\gamma$ ensures an expected search time $T(\gamma, x)$ of
\begin{align}
T(\gamma, x) \le \mu_2 + V(Q^t) = \mu_2 + \frac{1}{2}(\mu_1 + D(Q^t)), \label{eq:gamma}
\end{align}
by Theorem~\ref{treegame}. To see that inequality~(\ref{eq:gamma}) also holds for any point $x$ that lies in a block of $Q$, we add a leaf arc of infinitesimal length incident to some point on the block. Clearly point $x$ is reached at an earlier time that the leaf arc, so since~(\ref{eq:gamma}) holds for points on the leaf arc it must also hold for points in the block.

\subsection{Performance of $\gamma$}
Our main result of this section is on the performance of the search strategy $\gamma$. 
\begin{thm}\label{thm:gamma}
	Let $Q$ be a rooted network with bridge-ratio $r$. The maximum expected search time $T(\gamma)$ of the bridge-optimal strategy $\gamma$ on $Q$ satisfies
	\[
	T(\gamma) \le \left( \frac{2}{1+r^2} \right) V(Q).
	\]
	That is, the bridge-optimal strategy $\gamma$ is $2/(1+r^2)$-factor optimal.
\end{thm}
In order to prove the theorem, we need to give a lower bound on $V(Q)$ that is tight in the case that $Q$ is a tree. We will do that by describing a general Hider distribution (probability measure) on $Q$. We note that the lower bound on $V$ obtained in this subsection applies
as well to pathwise search games, as $V\leq V^{P}$, so the results obtained
here are stronger, although sometimes the same proofs work.

We can improve on the uniform Hider distribution by placing all the uniform
measure on a bridge at its forward end (away from the root). In
fact we can do better than this. Moving upwards in $Q_1$ (starting
at arcs at the root), we successively remove the uniform measure on each such
arc $a$ and suitably increase the total mass of the uniform measure on the subtree of $Q_1$
following $a$ so that the total measure of $Q_1$ does not change. In
particular, when $a$ is a leaf arc of $Q_1$, we place its measure on the
corresponding leaf node. We call the resulting distribution the \emph{pushed
	uniform distribution} and denote it by $u^{\ast}$. The measure on $Q_1$ is thus repeatedly pushed away from the root until it is all on the leaf nodes. Note that $u^{\ast}$ agrees
with the uniform distribution $u$ on the blocks of $Q$ and the rest is concentrated on the
leaf nodes of $Q_1$.

\begin{figure}[ht]
	\begin{center}
		\includegraphics[scale=0.5]{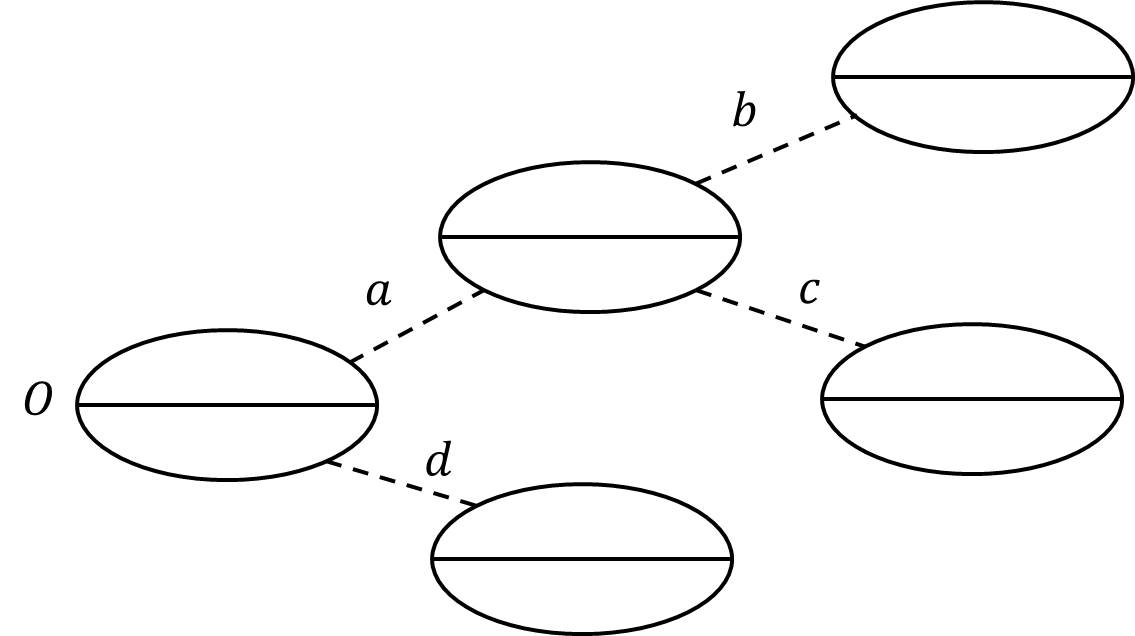}
		\caption{Decomposition of a particular network.}\label{fig:decomp}
	\end{center}
\end{figure}

For example on any circle-with-spike network $Q_{\alpha}^{CS},$ $u^{\ast
}\left(  B\right)  =1/3$, $u^{\ast}$ is uniform on its circle, and it is zero
on the open arc $AB$. Note that if $Q$ is already a tree, then $u^{\ast}$ is
the EBD distribution. Figure~\ref{fig:decomp} depicts a network all of whose arcs have unit length and whose
2-arc-connected components are copies of the three-arc network, given by two nodes joined by three unit length arcs. The network has 19 unit
length arcs, so each has a uniform measure of density $1/19$ under the uniform
measure $u$. The five three-arc networks keep this measure under $u^{\ast}$.
The measure $1/19$ on $a$ is split between $b$ and $c$, and eventually the
nodes at the forward ends of $b$ and $c$ each get measure $1/19+1/28=3/28$.
The measure $1/19$ on $d$ gets transferred to the forward end of $d$. Thus the
three leaf nodes of $Q^t$ get measure $3/28,$ $3/28,$ and $2/28$ which is
exactly $\mu_{1}/\mu=4/19$ times the EBD measure on $Q_1$. In general,
$u^{\ast}\left(  i\right)  =\left(  \mu(Q_1)/\mu\right)  \cdot e\left(
i\right)  ,$ where $e$ is the EBD measure on $Q^t$, for any leaf node $i$ of
$Q_1$.

\begin{lemma}
	\label{moduniformtheorm}For any network $Q$, let $\mu_{1}=\mu (Q_1)$. Then
	\begin{equation}
	V\geq \frac{\mu +\left( \mu _{1}/\mu \right) ~D(Q^t) }{2}%
	.  \label{moduniformestimate}
	\end{equation}
\end{lemma}

\begin{figure}[th]
	\begin{center}
		\includegraphics[scale=0.5]{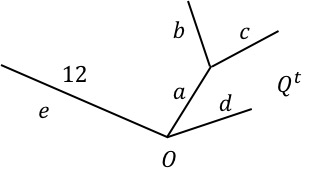}
	\end{center}
	\caption{The network $\tilde{Q}$.}%
	\label{fig:Qtilde}%
\end{figure}

\proof{Proof.} The Hider can ensure at least this expected capture time by
adopting the pushed uniform distribution $u^{\ast}$. Suppose we increase the
strategy space for the Searcher so that he can move freely within any given block and can thus search $Q_1$ as if it were
the tree $Q^t$. This is equivalent to changing the search space to the network
$\tilde{Q}$ consisting of the tree ${Q}^t$ and an additional
leaf arc $e$ at the root of length $\mu_2$. The network $\tilde{Q}$ corresponding to the network of Figure
\ref{fig:decomp} is shown in Figure \ref{fig:Qtilde}.

Note that any search $S$ of $Q$ induces a search $\tilde{S}$ of $\tilde{Q}$. At
time $t$, $\tilde{S}\left(  t\right)  $ will include a distance $\lambda\left(
S\left(  t\right)  \cap a\right)  $ along $a$, the amount of
$a$ that has been searched, as well as the identical portions of $a,b,c$
and $d.$ Considering $u^{\ast}$ as a distribution on $\tilde{Q}$ (which is
uniform on $a$ and EBD on $Q^t,$ we see that the subtrees
$a$ and $Q^t$ both have the maximum density of $1/\mu$. By
Lemma \ref{uniform} there is an optimal search that searches $a$ all at once, and by Lemma \ref{searchdensity}, $Q^t$ and
$a$ can be searched in either order, so assume $a$ is
searched first and then $Q^t$ is searched optimally. If the Hider is
located on $a$ the expected search time is half its length, $\mu_2/2$. So the expected search time is
\begin{align*}
& \frac{\mu_2}{\mu}\left(  \frac{\mu_2}{2}\right)  +\frac{\mu_1}{\mu
}\left(  \mu_2+V\left(  Q^t\right)  \right) \\
& =\frac{\mu_2}{\mu}\left(  \frac{\mu_2}{2}\right)  +\frac{\mu_1}{\mu
}\left(  \mu_2+\frac{\mu_1+D\left(  Q^t\right)  }{2}\right)
\text{ (by Theorem \ref{treegame})}\\
& =\frac{\mu+\left(  \mu_1/\mu\right)  ~D\left(  Q^t\right)  }{2}. 
\end{align*}
\endproof

To illustrate the construction in the proof of Lemma~\ref{moduniformtheorm}, we consider again the example of Figure~\ref{fig:Qtilde}, we have $\mu=19,\mu
_{1}=4$, and
\[
D(Q^t)  =3/8\cdot2+3/8\cdot2+2/28\cdot1=14/8,
\]
so the value satisfies
\[
V\left(  Q\right)  \geq\left(  19+\left(  4/19\right) \left(  14/8\right)  \right)  /2=184/19 = 9.6842.
\]
We note that we have equality in (\ref{moduniformestimate}) when $Q$ is a tree
($\mu_{1}=\mu$) or if $Q$ is $2$-arc-connected ($\mu_{1}=0$), by Theorem~\ref{treegame} and Theorem~\ref{thm:2arcconnected}.

We can now prove Theorem~\ref{thm:gamma}.

\proof{Proof of Theorem~\ref{thm:gamma}.}
Combining the lower bound in~(\ref{moduniformestimate}) with our estimate~(\ref{eq:gamma}) for the expected search time of the search strategy $\gamma$, we have
\begin{align*}
\frac{T(\gamma)}{V} &\le \frac{2\mu_2+(\mu_1+D(Q^t))}{\mu+(\mu_1/\mu)D(Q^t)} \\
&= \frac{(2-r)\mu+D(Q^t)}{\mu+rD(Q^t)} \mbox{ (where $r=\mu_1/\mu$)}\\
&\le \frac{(2-r)\mu+r\mu}{\mu+r^2\mu} \mbox{ (since $D \le \mu_1$)}\\
&=\frac{2}{1+r^2}. 
\end{align*}
\endproof

We have proved that the strategy $\gamma$ is $2/(1+r^2)$-factor optimal. This factor is decreasing in $r$: when $r=1$ (so that $Q$ is a tree), then $\gamma$ is optimal; when $r=0$ (so $Q$ is $2$-arc-connected), then $T(\gamma)=2V$ and $\gamma$ performs very badly.

\section{Summary and Conclusions}
\label{sec:conclusion}

In many or even most searches that are carried out to find missing persons,
lost airplanes or unexploded mines from earlier conflicts, the rate at which
the searched area can expand is restricted by available resources. In such
situations the ``expanding search'' model recently introduced by the authors
in \cite{AL13} can be applied. This model seeks the
randomized strategy for expanding the search region which minimizes the
expected time to find the ``target'' in the worst case, when the search region
has a network structure. Optimal strategies in this context were previously
known only for networks that are trees or are $2$-arc-connected. This paper
extends the classes of solvable networks to include the new class of
``circle-with-spike'' networks. However our main contribution is to give two
general classes of strategies which can be applied to {\em any} network
and have expected times to find the target that are within 20\% of the
optimal time. One class of strategies, which was introduced and analyzed in
Section~\ref{sec:approx}, is the block-optimal strategy $\beta $, which is optimal when the
network has no bridge (when the bridge-ratio $r$ is equal to $0$). This
strategy is also close to optimal when the ``height'' $\pi $ of the network is
close to $1$ (when one of the bridges is much longer than all the others
combined). The other class of strategies, introduced in Section~\ref{sec:treelike}, is called the bridge-optimal strategy 
$\gamma$, and is optimal when when the search network is a tree. It is close to
optimal when the bridge-ratio $r$ of the network is close to $1$. 

The paper enables the organizer of an expanding search effort to base his
search strategy simply on the bridge-ratio $r$ of the network to be
searched. This is the ratio $r=\mu _{1}/\mu $ of the total length of the
bridges of the network to the total length of all its arcs. Roughly
speaking, when $r$ is small the searcher can get closest to the optimal
expected search time by adopting the block-optimal strategy $\beta $ and
when $r$ is large (above 80\%) by adopting the bridge-optimal strategy $%
\gamma$. 

\subsection{Factor-approximate estimates in terms of the bridge-ratio $r$}

It is useful to combine the estimates on the value $V=V(Q)$
obtained by the block-optimal strategy~$\beta $ and the bridge-optimal
strategy~$\gamma$ in terms of the bridge ratio $r=\mu _{1}/\mu $ of the
network $Q$. Although we still get the factor $\left( 1+\sqrt{2}\right) /2$
for {\em all} networks, we can get better estimates for values of the
bridge ratio $r$ near $0$ or $1$, that is for block-like or bridge-like
networks. (Of course even better estimates can be obtained if we make the
more difficult calculation of the height $\pi$ of the network.)

First note that, as explained in the proof of Theorem~\ref{thm:approx}, $\pi/\mu \leq r$. Since the function $f(r)=(1+r)/(1+r^2)$ is increasing on the interval $\left[ 0,\sqrt{2}-1\right] \simeq \left[ 0,0.414\,21\right] $ with a maximum on $[0,1]$
of $(1+\sqrt{2})/2\simeq 1.2071$
obtained at $r=\sqrt{2}-1$, it follows that 
\begin{eqnarray*}
	f\left( \frac{\pi }{\mu }\right)  &\leq &f\left( r\right) \text{ for }r\in %
	\left[ 0,\sqrt{2}-1\right] \text{ and} \\
	f\left( \frac{\pi }{\mu }\right)  &\leq &\left( 1+\sqrt{2}\right) /2\text{
		for all }r\in \left[ 0,1\right] .
\end{eqnarray*}%
By Theorem~\ref{thm:approx}, we know that by adopting the block-optimal strategy $\beta $
the Searcher ensures that the ratio of actual expected time to optimal
expected time is bounded by 
\[
\frac{T}{V}=\frac{T\left( \beta \right) }{V}\leq f\left( \frac{\pi }{\mu }%
\right) .
\]%
Similarly from Theorem~\ref{thm:gamma}, we know that by adopting the bridge-optimal
strategy $\gamma $ the Searcher ensures that 
\[
\frac{T}{V}=\frac{T\left( \gamma \right) }{V}\leq \frac{2}{1+r^{2}}.
\]%
Observe that $g(r)=2/(1+r^{2})$ is decreasing in $r
$ and is equal to $(1+\sqrt{2})/2\simeq 1.207$ for $r=\sqrt{\frac{3-\sqrt{2}}{1+\sqrt{2}}}\simeq 0.81047$.

So it follows that by taking strategy $\alpha $ to be the better of the
block-optimal and bridge-optimal strategies $\beta $ and $\gamma ,$ we can
ensure that%
\[
\frac{T\left( \alpha \right) }{V}\leq \left\{ 
\begin{array}{ccc}
f\left( r\right) =\frac{1+r}{1+r^{2}} & \text{ for} & 0\leq r\leq \sqrt{2}%
-1\simeq 0.414\,21, \\ 
f\left( \sqrt{2}-1\right) =\frac{1+\sqrt{2}}{2}\simeq 1.\,\allowbreak 207\,1
& \text{for} & \sqrt{2}-1\leq r\leq \sqrt{\frac{3-\sqrt{2}}{1+\sqrt{2}}}%
\simeq 0.810\,4, \\ 
g\left( r\right) =\frac{2}{1+r^{2}} & \text{for} & \sqrt{\frac{3-\sqrt{2}}{1+%
		\sqrt{2}}}\leq r\leq 1.%
\end{array}%
\right. 
\]%
The first two estimates are obtained by adopting $\alpha =\beta $ and the
last one by adopting $\alpha =\gamma .$ The upper bounds on $T/V,$ as a
function of the bridge ratio $r,$ are shown in Figure 5. 

\begin{figure}[th]
	\begin{center}
		\includegraphics[scale=0.45]{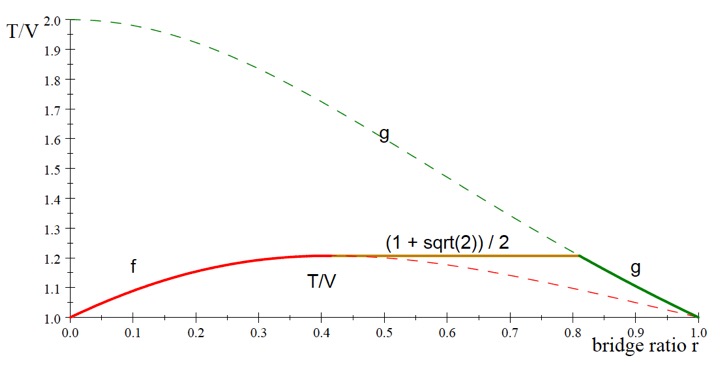}
	\end{center}
	\caption{$T/V$ bounded above by thick line, as a function of the bridge ratio $r$.}%
	\label{fig:graph}%
\end{figure}

We can see from Figure~\ref{fig:graph} that by
adopting one of only the two expanding search strategies $\beta $ and $%
\gamma $ introduced in this paper, a searcher with an expanding search to
plan on a known network can be sure of getting within $20\%$ of the optimal
expected time when the network has a bridge-ratio $r$ within $0.41\leq r\leq
0.81$ and can be even closer to optimal for $r$ outside that interval. By
considering the height $\pi $ of the network, he might be able to ensure
being even closer to optimal, but the calculation of $r$ is very easy. It is
worth noting that prior to this paper there were no known general strategies
for expanding search on arbitrary networks. 

\subsection{Conclusions}

From the results of this paper, it seems that there two directions for
further investigations in the study of expanding search on networks. One is
to identify more classes of networks where optimal strategies can be found.
The other direction is to see if some strategy classes can be found which
reduce the factor-approximate constant of $1.2071$ for the ratio of $T/V.$
Another aspect of real-life search that might be added to the analysis is a
``give-up'' option, that is, when should the search be stopped. This is
particularly relevant to the Bayesian problem where there is an {\em a priori}
distribution for the target, as in the search for the Malaysian Airlines
plane.

\section*{Acknowledgements}
Steve Alpern wishes to acknowledge support from the Air Force Office of Scientific Research [Grant FA9550-14-1-0049].


\end{document}